\documentclass{amsart}
\usepackage{CJKutf8}
\usepackage{amsmath,amsfonts,amssymb,amscd,bm,latexsym}
\usepackage[colorlinks=true]{hyperref}
\usepackage{tikz}
\usepackage{pgflibraryarrows}
\usepackage{pgflibrarysnakes}
\usepackage{pgflibraryarrows}
\usepackage{pgflibrarysnakes}
\hypersetup{urlcolor=blue, citecolor=blue}
\numberwithin{equation}{section} 

\def\<{\langle}             \def\>{\rangle}

\newtheorem{thm}{Theorem}[section]

\newtheorem{lem}[thm]{Lemma}

\theoremstyle{definition}

\newcommand{\beeq}{\begin{equation}}\newcommand{\eneq}{\end{equation}}
\newcommand{\al}{\alpha}    \newcommand{\be}{\beta}
\newcommand{\de}{\delta}    
\newcommand{\ep}{\varepsilon}
  
    \newcommand{\la}{\lambda}

\newcommand{\om}{\omega}    \newcommand{\Om}{\Omega}
\newcommand{\ga}{\gamma}    
\newcommand{\R}{\mathbb{R}}

\newenvironment{prf}{\noindent {\bf Proof.} }{\endprf\par}
\def \endprf{\hfill  {\vrule height6pt width6pt depth0pt}\medskip}
\newcommand{\pa}{\partial}
\newcommand{\les}{{\lesssim}}
\newcommand{\gt}{{\gtrsim}}
\newcommand{\supp}{\,\mathop{\!\mathrm{supp}}}

\numberwithin{equation}{section}
\title[Wave equations in Schwarzschild space]
      {Blow-up for semilinear wave equations with damping and potential in high dimensional Schwarzschild spacetime}
\author{Mengliang Liu}
\address{Department of Mathematics\\                	
Guangxi University\\                Nanning 530000, P. R. China}
\email{l15077158750@163.com}
\author{Mengyun Liu$^{*}$}\thanks{* Corresponding author}
\address{Department of Mathematics\\Zhejiang Sci-Tech University\\Hangzhou 310018, P. R. China }
\email{mengyunliu@zstu.edu.cn}

\date{\today}
\dedicatory{} \commby{}

\thanks{The second author was supported by NSFC 12101558 and NSF of Zhejiang province LQ22A010016.
}

\begin{document}
\begin{abstract}
In this work, we study the blow up results 
to power-type semilinear wave equation in the high dimensional Schwarzschild spacetime, with damping and potential terms. We can obtain  the upper bound estimates of lifespan without the assumption that the support of the initial date should be far away from the black hole.
\end{abstract}


\maketitle

\section{Introduction}
In this paper, we are interested in the blow up of solutions for the following semilinear wave equation with small initial data, posed on the high dimensional Schwarzschild black hole 
\begin{align}
\begin{cases}
\label{h1}
\square_{g_{S}}u+D(r)u_{t}+V(r)u=G(u, u_{t}), \ (0, T) \times \Omega,\\
u(x,0)=\ep u_{0}(x),\ \ \ \ u_{t}(x,0)=\ep u_{1}(x),\ x\in \Om,
\end{cases}
\end{align}
where $g_{S}$ denotes the Schwarzschild metric and the D'Alembert operator $\Box_{g_S}$ associated with the Schwarzschild metric has the form
\beeq
\square_{g_{S}}=\frac{1}{F(r)}\Big(\pa^{2}_{t}-\frac{F(r)}{r^{n-1}}\pa_{r}\big(r^{n-1}F(r)\big)\pa_{r}-\frac{F(r)}{r^{2}}\Delta_{\mathbb{S}^{n-1}}\Big)
\eneq
with $F(r)=1-2M/r$, $\Delta_{\mathbb{S}^{n-1}}$ is the standard Laplace-Beltrami operator on $\mathbb{S}^{n-1}$. Here, $\Omega= \{(r,\omega) : r>2M,\  \omega \in \mathbb{S}^{n-1}\}= (2M, \infty) \times \mathbb{S}^{n-1}$, $n\geq3$. 
For the initial data, we assume they are nonnegative, nontrivial, and compactly supported
\beeq
\label{d1}
 \supp\  u_{0}(x), u_{1}(x) \subset \{2M+R_{1}\leq r \leq 2M+R_{2}\}\times \mathbb{S}^{n-1},
\eneq
where $R_{1}<R_{2}$ are some positive constants.
Concerning the damping and potential terms, we assume $D(r), V(r)$ takes the form 
\beeq\label{khbd00}
0\leq D\leq \mu_1r^{-\be}\ ,\ 0\leq V\leq \mu_2 r^{-\ga}\ , \ \be>1\ , \ga>2\ ,
\eneq 
for some $\mu_1, \mu_2\geq 0$. 
For the nonlinear terms $G(u.u_t)$, we mainly consider the two typical terms $|u|^{p}$, $|u_t|^{p}$.

The first solution of the Einstein equation was obtained by Schwarzschild in 1916, with the metric given by 
\beeq
\label{metric}
g_{S}=(1-\frac{2M}{r})dt^{2}-(1-\frac{2M}{r})^{-1}dr^{2}-r^{2}d\omega^{2}
\eneq
in standard spherical coordinates. Here, the constant $M>0$ denotes the mass of the black hole and $d\om^{2}$ is the standard metric on the unit sphere $\mathbb{S}^{2}$. 
The metric has two type singularities, the singularity at $r=2M$ is only a coordinate singularity (event horizon). In fact, if we use tortoise coordinate
\beeq
\label{wg}
s=r+2M\ln(r-2M)\ ,
\eneq
the metric becomes
$$g_S=F\Big(dt^2-ds^2\Big)-r^2d\om^2\ .$$
But $r=0$ is the first example of a black hole singularity, curvature will blow up when $r$ tends to zero. 

Before presenting our main results, let us first give a brief review of the history, in a broader context.

\noindent
Notations. We will use $A\les B$ $(A\gtrsim B)$ to stand for $A\leq CB$ $(A\geq CB)$ where the constant C may change from line to line. We also use $A\simeq B$ to stand for $c_1B\leq A\leq c_2B$ for some positive constants $c_1\leq c_2$.
\\

{\setlength{\parindent}{1em}\bf (\uppercase\expandafter{\romannumeral1}) Minkowski spacetime, $G(u,u_{t})=|u|^{p}$}\\

When $g_S=m_{\al\be}=diag(1, -1, -1, \cdot\cdot\cdot, -1)$ and $\Om=\R^n$, \eqref{h1} is related to the Strauss conjecture
\beeq
\label{srauss1}
\pa^{2}_t u-\Delta u+Du_t+Vu=|u|^p\ , (t, x)\in (0, T)\times\R^n\ .\eneq
It was conjectured that \eqref{srauss1} admits a critical power
$p_S(n)$ when $D=V=0$, 
where $p_S(n)$ is the positive root of the quadratic equation
\beeq
\ga(p,n) :=2+(n+1)p-(n-1)p^{2}=0\ .
\eneq
This conjecture has been essentially proved,
see, e.g., \cite{F, R, Y-Z, H-C, V-H-C, N-L-Y} for global results and \cite{F, R-T, J, T, B-Q-Z} for blow up results (including the critical case $p=p_{S}(n)$). In recent years, there are many works concerning the influence of damping and potential terms to \eqref{srauss1}. For the space dependent damping terms, we refer \cite{L-L-W-T} for the history. In particular, for the typical damping $D=1/(1+|x|)^{\be}$ and potential $V=1/(1+|x|)^{\ga}$, Liu-Lai-Tu-Wang \cite{L-L-W-T} showed the finite time blow results under the Strauss power $p<p_S(n)$ when $\be>1$ and $\ga>2$.

Besides, there are many works studied the time dependent damping terms see, e.g., \cite{T-H-T-O, P-K, K-N, J-W, L-H} and references therein. 
\\

{\setlength{\parindent}{1em}\bf (\uppercase\expandafter{\romannumeral2}) Schwarzschild spacetime, $G(u,u_{t})=|u|^{p}$}\\

When there are no damping and potential, the analogous of Strauss conjecture on 3-dim Schwarzschild spacetime has been proved. Lindblad-Metcalfe-Sogge-Tohaneanu-Wang \cite{L-W} established global existence when $p>1+\sqrt2$. For the long time existence, Wang \cite{W} established the lower bound estimates of the lifespan for many spacetimes including the nontrapping exterior domain, nontrapping asymptotically Euclidean space and Schwarzschild spacetime. For the blow up part, 
Catania-Georgiev \cite{D-V-G} proved blow-up result in some weak sense when $1<p<1+\sqrt2$. Lin-Lai-Ming \cite{L-N} obtained the blow up results when $1<p\leq 2$. Recently,  Lai-Zhou \cite{N-Y} obtained the blow up results for $2<p\leq1+\sqrt2$ and
\begin{align*}
T_\ep\ \les\ 
\begin{cases}
\ep^{-\frac{p(p-1)}{1+2p-p^{2}}},\ &2\leq p<1+\sqrt2\ ,\\
\exp \big(\ep^{-(2+\sqrt2))}\big),\ &p=1+\sqrt2\ ,
\end{cases}
\end{align*}
with the support of the initial data can be close the event horizon. Here, $T_\ep$ denotes the lifespan and $\ep$ is the size of the initial data. Besides, Liu-Wang \cite{MLL} showed the blow up phenomenon to \eqref{h1} on a large class of asymptotically flat space including Schwarzschild and Kerr black hole backgrounds under the Strauss exponent $1<p<1+\sqrt{2}$. \\
 
{\setlength{\parindent}{1em}\bf (\uppercase\expandafter{\romannumeral3}) Minkowski spacetime, $G(u,u_{t})=|u_{t}|^{p}$}\\

When $g_S=m_{\al\be}=diag(1, -1, -1, \cdot\cdot\cdot, -1)$, $\Om=\R^n$ and $D=V=0$, \eqref{h1} is related to the Glassey conjecture
$$\pa^{2}_t u-\Delta u=|u_t|^p\ , (t, x)\in (0, T)\times \R^n\ ,$$
which was proposed in \cite{R-T-T}. The problem admits a critical power
\begin{align*}
p_{G}(n)=
\begin{cases}
\frac{n+1}{n-1},\ \ n\geqslant2\ ,\\
\infty,\ \ \ \ \ n=1\ ,
\end{cases}
\end{align*}
and this conjecture has been essentially proved, see, e.g.,\cite{LML, J-S, R-A, T-C, Y,  Z-W} for the history. Concerning the typical damping $D=(1+|x|)^{-\be}$, Lai-Tu \cite{L-T} obtained the upper bound estimate of the lifespan with $\be>2$ and $n\geq 2$ 
\begin{align}
\label{502}
T\ \les\ 
\begin{cases}
\ep^{-(\frac{1}{p-1}-\frac{n-1}{2})^{-1}}\ ,\ &1<p<p_{G}(n)\ ,\\
\exp (\ep^{-(p-1)})\ ,\ &p=p_{G}(n)\ .
\end{cases}
\end{align}

{\setlength{\parindent}{1em}\bf (\uppercase\expandafter{\romannumeral4}) Schwarzschild spacetime, $G(u,u_{t})=|u_{t}|^{p}$}\\

When $n=3$ and $D=V=0$, Lai-Zhou \cite{L-Z} obtained the solution of \eqref{h1} will blow up in finite time with $1<p\leqslant2$, they also established the upper bound estimate of the lifespan
\begin{align}
T \ \les\ 
\begin{cases}
\ep^{-\frac{p-1}{2-p}},\ \ \ 1<p<2\ ,\\
\exp \big(\ep^{-1}),\ \ \ \ \ p=2\ .
\end{cases}
\end{align}

 There is not much works concerning the dampings and potentials posed on black hole spacetime. In this work, we show that the solution of \eqref{h1} with damping and potential terms \eqref{khbd00} will blow up in finite time for any $n\geq 3$.
 

Note that the assumption \eqref{d1} on the initial data is equivalent to 
\beeq
\label{supp01}
\supp( u_0(r(s)), u_1(r(s)))\subset\ \{s\in \R;\  |s|\leq R\}\times\mathbb{S}^{n-1}\ ,
\eneq
for some constant $R>0$. As usual, to show the finite time blow up results of wave equations, we need the solution to satisfy finite speed of propagation. For that purpose, we shall make a hypothesis.
\\
{\bf Hypothesis:} We assume the weak solution of \eqref{h1} satisfying 
\beeq
\label{fps}
\supp u\ \subset \{|s|\leq t+R\}\times \mathbb{S}^{n-1}\ ,
\eneq
for any $t\in (0, T_\ep)$.

\begin{thm}
\label{thm1}
Let $n\geq3$, $D, V\in C(2M, \infty)$, $G(u,u_{t})=|u|^{p}$. Consider the Cauchy problem \eqref{h1} with \eqref{d1} and \eqref{khbd00}. Then for any $1<p< p_{S}(n)$, any weak solutions of \eqref{h1} that satisfying \eqref{fps} will blow up in finite time. 
In addition, there exist positive constants $C$, $\ep_0$ such that the lifespan $T_{\ep}$ satisfies
\begin{align}
\label{d4}
T_{\ep}\leq
\begin{cases}
C\ep^{-\frac{2(p-1)}{n+1-(n-1)p}},\ \ &1<p\leq \frac{n}{n-1},\\
C\ep^{\frac{-2p(p-1)}{\ga(p,n)}},\ \ &\frac{n}{n-1}<p<p_{S}(n),
\end{cases}
\end{align}
for any $0<\ep\leq \ep_0$. In addition, when $V=0$, we have
\beeq
T_\ep \ \leq\ \exp(C\ep^{-p(p-1)}),\ \ p=p_{S}(n)\ .
\eneq
\end{thm}

\begin{thm}
\label{thm2}
Let $n\geq3$, $D, V\in C(2M, \infty)$, $G(u,u_{t})=|u_{t}|^{p}$. 
 Consider the Cauchy problem \eqref{h1} with \eqref{khbd00}. Suppose the initial data satisfies \eqref{supp01} and 
 \beeq
 \label{chuzhi2}
\int_{\R}\phi_{\la_0}\bigg(u_1 r(s)^{(n-1)/2}+\big(\la_0 DF+\la_{0}^{2}-Q\big)u_0r(s)^{(n-1)/2}\bigg)ds >0\ , 
 \eneq
where $Q$ is in \eqref{R1} and $\phi_{\la_0}$ is in Lemma \ref{le4} with some fixed $\la_0>4/(M(p-1)p)$. 
  Then for any $1<p\leq p_{G}(n)$, any weak radial solutions of the \eqref{h1} that satisfying \eqref{fps} will blow up in finite time. 
In addition, there exist positive constants $C$, $\ep_0$ such that the lifespan $T_{\ep}$ satisfies
\begin{align}\label{c9}
T_{\ep}\leq
\begin{cases}
C\ep^{-(\frac{1}{p-1}-\frac{n-1}{2})^{-1}},\ \ &1<p<p_{G}(n),\\
\exp \big(C\ep^{-(p-1)}\big),\ \ &p=p_{G}(n).\\
\end{cases}
\end{align}
for any $0<\ep \leq \ep_{0}$.
\end{thm}


\section{PRELIMINARY}
In this section, we collect some Lemmas we will use later. For the strategy of the proof, we basically follow the test function methods (see, e.g., \cite{Q-Z}, \cite{B-Q-Z}, \cite{I-M-S-M-W-K}). For simplicity, we define 
\beeq
\Delta_{g_{S}}=\frac{1}{r^{n-1}}\pa_{r}\Big((r^{n-1}F)\pa_{r}\Big)+
\frac{1}{r^{2}}\Delta_{\mathbb{S}^{n-1}}\ , 
\eneq 
then the linear equation of \eqref{h1} takes the form 
$$\frac{1}{F}\pa^{2}_t u-\Delta_{g_{S}}u+D\pa_t u+Vu=0\ .$$
The key ingredient of the test function method is to construct the special positive solutions of the form $\Psi=e^{-\la t}\psi_{\la}$ to the linear dual problem 
$$\frac{1}{F}\pa^{2}_t \Psi-\Delta_{g_{S}}\Psi-D\pa_t \Psi+V \Psi =0\ ,$$
with desired asymptotic behavior. In turn, we just have to construct solutions to the corresponding elliptic problems
\beeq
 \Delta_{g_{S}}\psi_{\la}=\frac{1}{F}\la^2\psi_{\la}+\la D\psi_{\la}+V \psi_{\la} \ .
 \eneq
To describe the behavior of $\psi_{\la}$, we need to introduce the variable 
\beeq
s=r+2M\ln(r-2M)\ ,
\eneq
which is known as tortoise coordinate.  Furthermore, it is easy to see that 
\begin{align}
\label{srelation}
\begin{cases}
s\leq 4M+e\Longleftrightarrow2M\leqslant r\leqslant2M+e \ \Longleftrightarrow r-2M\simeq e^{\frac{s}{2M}}\ .\\
s\geq 4M+e\Longleftrightarrow r \geq 2M+e \ \Longleftrightarrow s\simeq r\ .
\end{cases}
\end{align}

\begin{lem}
\label{le6}
Let $n\geq3$. Suppose $V\in C(2M, \infty)$ and satisfying \eqref{khbd00}. Then the elliptic equation
\begin{align}
\label{110}
\Delta_{g_{S}}\psi_0(x)=V\psi_{0}(x)\ , \ x\in \Om\ , 
\end{align}
exists a $C^{2}$ radial solution satisfying
\beeq
\psi_{0}(r) \simeq1, \ \ \pa_r\phi_0(r)\ \les \ \frac{1}{r}\ , \ r> 2M\ .
\eneq

\end{lem}
\begin{prf}
Let $\psi_{0}(x)=\psi_0(|x|)=\psi_{0}(r)$, then  \eqref{110} becomes an ODE  
$$\frac{1}{r^{n-1}}\pa_{r}\Big((r^{n-1}F)\pa_{r}\psi_{0}(r)\Big)=V(r)\psi_{0}(r).$$
We consider the following initial value problem
\beeq
\label{617}
\begin{cases}
\psi_{0}''+(\frac{(n-1)F}{r}+\pa_{r}F)\psi_{0}'-V\psi_{0}=0,\ r>2M\ , \\
\psi_{0}(2M)=1, \psi_{0}'(2M)=0 .
\end{cases}
\eneq
Then by the theory of ordinary differential equation, there is a unique solution $\psi_{0}\in C^2(2M, \infty)$. 
Recall that
$$
\label{2-d1}
\pa_r(r^{n-1}F\pa_r\psi_{0})=r^{n-1}V\psi_{0}\geq0 , 
$$
we see that $\phi_0\geq 0$ and $\pa_r\phi_0\geq 0$ for any $r\geq 2M$. 
By integrating it from $2M$ to $r$, we get
$$
r^{n-1}\pa_r\psi_{0} F=\int^{r}_{2M}\tau^{n-1}V\psi_{0} d\tau \ \les \ \psi_{0}\int^{r}_{2M}(1+\tau)^{n-1-\be}d\tau\ .
$$
Thus, as $\be>2$, if $n\geq 3$ and $r\ge 2M$, we have
\beeq
\label{psi1}
\pa_r\psi_{0}\ \les \
 \psi_{0} r^{1-n}\int^{r}_{2M}(1+\tau)^{n-1-\be}d\tau\ \les\ 
\psi_{0} r^{-1-\de}\ ,
\eneq
for some $\de>0$. 
By Gronwall's inequality, we obtain
$\psi_{0}(r) \les  1$, for all $r\geq 2M$,
which yields $$\psi_{0}\simeq 1\ .$$
By \eqref{psi1}, we get that $\pa_r\psi_0\les 1/r$. 
\end{prf}



\begin{lem}
\label{le5}
Let $n\geq3$. Suppose $D, V\in C(2M, \infty)$ and satisfying \eqref{khbd00}. Then for any $0<\la\in\R$, the eigenvalue problem 
\begin{align}
\label{b2}
\Delta_{g_{S}}\psi_{\la}=\frac{\la^{2}}{F}\psi_{\la}+D\la\psi_{\la}+V\psi_{\la}, \ x\in \Om\ , 
\end{align}
admits a positive radial solution $\psi_{\la}$ with
\begin{align}
\label{b4}
\psi_{\la}(r)\simeq r^{-\frac{n-1}{2}}e^{\la s}\ , \ |\pa_{r}\psi_{\la}(r)|\les \frac{\la+1}{r^{\frac{n-1}{2}}F}e^{\la s}\ , \ r>2M\ .
\end{align}
In addition, when $V=0$, we have 
\beeq
\ |\pa_{r}\psi_{\la}(r)|\ \les\  \frac{\la \psi_{\la}}{F}\ , \ r>2M\ .
\eneq
\end{lem}
\begin{prf}
By the tortoise coordinate we see that $F(r)\pa_r =\pa_s$, then \eqref{b2} becomes
$$\frac{1}{F}\pa^{2}_{s}\psi_\la+\frac{(n-1)}{r(s)}\pa_{s}\psi_\la+\frac{1}{\big(r(s)\big)^{2}}\Delta_{\mathbb{S}^{n-1}}\psi_\la=\frac{\la^{2}}{F}\psi_{\la}+D\la\psi_{\la}+V\psi_{\la}\ .$$
Let $\psi_{\la}(x)=\psi_\la(r)$ and $\phi_\la(s)=\psi_{\la}(r(s))r(s)^{\frac{n-1}{2}}$, then $\phi_\la(s)$ satisfying the following ODE
$$\pa_{s}^{2}\phi_{\la}=Q(s)\phi_{\la}+\la D(r(s))F(r(s))\phi_{\la}+\la^{2}\phi_{\la}\ ,$$
where $Q(s)$ is in \eqref{R1}. 
By Lemma \ref{le4}, we know that $\phi_{\la}(s)\simeq e^{\la s}$. Thus we obtain that
$$\psi_\la(r)\simeq \ r^{-\frac{n-1}{2}}e^{\la s}\ .$$
For the derivative of $\psi_\la$, we know that
\begin{align}
\nonumber
\pa_{r}\psi_{\la}
=&\frac{\la^{2}}{r^{n-1}F}\int^{r}_{2M}\frac{\tau^{n}}{\tau-2M}\psi_{\la}d\tau+\frac{\la}{r^{n-1}F}\int_{2M}^{r}\tau^{n-1}D\psi_{\la}d\tau\\\label{q1}&+\frac{1}{r^{n-1}F}\int^{r}_{2M}V\psi_{\la}\tau^{n-1}d\tau\geq 0\ ,
\end{align}
for any $r>2M$. In addition, let $h=\tau+2M\ln(\tau-2M)$, then $d\tau=\frac{\tau-2M}{\tau}d h$. Note that when $r\in(2M, 2M+e]$ we have $s\leq 4M+e$ and $\tau\in (2M, r)\subset (2M, 2M+e)$, then we have
\begin{align*}
&\int^{r}_{2M}\frac{\tau^{n}}{\tau-2M}\psi_{\la}d\tau \les\  \int^{s}_{-\infty}\tau^{n-1}\psi_\la dh \les\ \frac{e^{\la s}}{\la}\ , \\
&\int^{r}_{2M}V\psi_{\la}\tau^{n-1}d\tau\leq \int_{2M}^{r}\tau^{n-1}D\psi_{\la}d\tau
 \les \ \psi_\la \int^{s}_{-\infty}e^{h/(2M)}dh \les e^{\la s}\ .
\end{align*}
 Back to \eqref{q1}, we get that 
 $$\pa_r\psi_{\la}(r)\les \frac{\la+1}{r^{(n-1)}F}e^{\la s}\ , r\in(2M, 2M+e]\ .$$
For the part $r>2M+e$, we have $s> 4M+e$, then 
\begin{align*}
\int^{r}_{2M}\frac{\tau^{n}}{\tau-2M}\psi_{\la}d\tau = &\int^{4M+e}_{-\infty}\tau^{n-1}\psi_\la dh+ \int_{4M+e}^{s}\tau^{n-1}\psi_\la dh\\
\les &\ \frac{e^{\la (4M+e)}}{\la}+s^{(n-1)/2}\frac{e^{\la s}}{\la}\ .
\end{align*}
Notice that $D(\tau(h))\simeq 1$ when $h\leq 4M+e$ and $D\les h^{-1}$ when $h\geq 4M+e$, then we have 
\begin{align*}
\int_{2M}^{r}\tau^{n-1}D\psi_{\la}d\tau=&\int^{4M+e}_{-\infty}D\tau^{n-2}(\tau-2M)\psi_\la dh+ \int_{4M+e}^{s}D\tau^{n-2}(\tau-2M)\psi_\la dh\\
\les&\ \psi_\la\Big(e^{\frac{4M+e}{2M}}+s^{n-1}\Big) \les \ s^{n-1}\psi_\la\ \ .
\end{align*}
We finish the proof of \eqref{b4}.
\end{prf}

As usual, to show the finite time blow up results when $p=p_S(n)$ and $V=0$, we need to construct a extra test functions
\beeq\label{bade}
b_{a}(t,r)=\int_{0}^{1}e^{-\la t}\psi_{\la}\la^{a-1}d\la\ , 0<a\in \R\ ,
\eneq
where $\psi_{\la}$ is in Lemma \ref{le5}\ . Then it is easy to check that $b_a$ satisfying 
\beeq
\frac{\pa}{\pa t}b_{a}=-b_{a+1}\ , \ 
\frac{1}{r^{n-1}}\pa_{r}(r^{n-1}F\pa_{r}b_{a})=\frac{1}{F}\pa_{t}^{2}b_{a}-D\pa_{t}b_{a}.
\eneq
Furthermore, we also need the lower bound as well as the upper bound estimates of $b_a$. 

\begin{lem}\label{leba}
Let $n\geq3$ and $0<a\in \R$. Then for $4M+e\leq s\leq t+R$, we have
\begin{align*}
&|\pa_{r}b_{a}(t, r)|\ \ \les \ b_{a+1}(t, r)\ , \ b_{a}(t, r)\ \simeq (t+R)^{-a}\ , \ 0<a<(n-1)/2\ , \\
&b_{a+1}(t, r)\les (t+R)^{-\frac{n-1}{2}}(t+R+1-s)^{-(a+1)}\ , \ a+1>(n-1)/2\ .
\end{align*}
\end{lem}
\begin{prf}
When $s\geq 4M+e$ we have that $r\geq 2M+e$, then we have $e/(2M+e)\leq F(r)\leq 1$, thus
\begin{align} \nonumber
|\pa_{r}b_{a}|\les&\int_{0}^{1}e^{-\la t}|\pa_{r}\psi_{\la}|\la^{a-1}d\la\les\int_{0}^{1}e^{-\la t}\psi_{\la}F^{-1}\la^{a}d\la \\
\les&\int_{0}^{1}e^{-\la t}\psi_{\la}\la^{a}d\la=b_{a+1}\ .
\end{align}

For the lower bound of $b_a$, since $\pa_r\psi_\la\geq 0$, we have that $\psi_\la(r)\geq \psi_\la(2M+e)$ for any $r\geq 2M+e$, then we get that
\begin{align*}
b_{a}\gtrsim&\int_{\frac{1}{2(t+R)}}^{\frac{1}{t+R}}\psi_{\la}(2M+e)e^{-\la(t+R)}\la^{a-1}d\la\\
\gtrsim&(t+R)^{-a}\int_{\frac{1}{2}}^{1}e^{-\kappa}\kappa^{a-1}d\kappa\\
\gtrsim&(t+R)^{-a}\ .
\end{align*}
To get the upper bound of $b_a$, we divide the proof into two parts: $4M+e\leq s\leq (t+R)/2$ and $(t+R)/2\leq s\leq t+R$. For the former case, we have 
\begin{align}
\label{b_a}
b_{a}\les \int_{0}^{1}e^{-\la (t+R)}e^{\la\frac{t+R}{2}}\la^{a-1}d\la \les(t+R)^{-a}\int_{0}^{\infty}e^{-\kappa}\kappa^{a-1}d\kappa\ \les(t+R)^{-a}.
\end{align}
If $(t+R)/2\leq s\leq t+R$, we have that
\begin{align*}
b_{a}\les& \int_{0}^{1}e^{-\la t}e^{\la(t+R)}(t+R)^{-\frac{n-1}{2}}\la^{a-1}d\la\\
\les&(t+R)^{-\frac{n-1}{2}}\int_{0}^{1}\la^{a-1}d\la\  \les(t+R)^{-\frac{n-1}{2}}.
\end{align*}
For the upper bound of $b_{a+1}$, by \eqref{b_a}, we have 
$$b_{a+1}\ \les\ (t+R)^{-a-1}\ , 4M+e\leq s\leq \frac{t+R}{2}\ .$$
If $(t+R)/2\leq s\leq t+R$, then for any $a+1>\frac{n-1}{2}$, we have
\begin{align*}
b_{a+1}\les&\int_{0}^{1}e^{-\la(t+R+1-s)}(\frac{t+R}{2})^{-\frac{n-1}{2}}\la^{a-\frac{n-1}{2}+\frac{n-1}{2}}d\la\\
\les&(t+R)^{-\frac{n-1}{2}}(t+R+1-s)^{-(a+1-\frac{n-1}{2})}\int_{0}^{\infty}e^{-\kappa}\kappa^{a-\frac{n-1}{2}}d\kappa\\
\les&(t+R)^{-\frac{n-1}{2}}(t+R+1-s)^{-(a+1-\frac{n-1}{2})}\ .
\end{align*}
\end{prf}

At the end of this section, we list an ODE lemma that we will use later. 
\begin{lem}[{\bfseries Lemma 3.10 in \cite{I-M-S-M-W-K}}] \label{le11}
Let $2<t_{0}<T $, $0\leq\phi\in C^{1}([t_{0},T))$. Assume that
\begin{align*}
\begin{cases}
\delta\leq K_{1}t\phi'(t),\  t\in(t_{0},T)\\
\phi(t)^{p_{1}}\leq K_{2}t(\ln t)^{p_{2}-1}\phi'(t), \  t\in(t_{0},T)
\end{cases}
\end{align*}
with $\delta$, $K_{1}$, $K_{2}>0$ and $p_{1}$, $p_{2}>1$.  If $p_{2}<p_{1}+1$, then there exists positive constant $\delta_{0}$ and $K_{3}$ (independent of $\delta$) such that
\begin{align*}
T\leq\exp\big(K_{3}\delta^{-\frac{p_{1}-1}{p_{1}-p_{2}+1}}\big)\nonumber
\end{align*} 
for all $\delta \in (0,\delta_{0})$.
\end{lem}

\section{Proof of Theorems \ref{thm1}}\label{proof}
In this section, we give the proof of Theorems \ref{thm1}.

\subsection{$1<p\leq \frac{n}{n-1}$}


We first introduce a smooth cut-off function $\eta(t)\in C^{\infty} ([0,\infty))$ such that
$$\eta(t)=1\ , \ t\leq 1/3\ , \ \eta(t)=0\ , \ t\geq 1\ , \ \eta'(t)\leq 0\ , t\geq 0\ .$$
Then for $T\in [1,T_{\ep})$, we set $\eta_{T}(t)=\eta(t/T)$. Let $\la=\la_0>\frac{4}{Mp(p-1)}$ to be taken. By multiplying the equation in \eqref{h1} with $\eta_{T}^{2p'}(t)e^{-\la_0 t}\psi_{\la_0}r^{n-1}$, where $\psi_{\la_0}$ is in Lemma \ref{le5} and making integration by parts, then we have
\begin{align} \nonumber
&C_{1}\ep+\int_{0}^{T}\int_{2M}^{\infty}\int_{\mathbb{S}^{n-1}}|u|^{p}\eta_{T}^{2p'}\psi_{\la_0}e^{-\la_0 t}r^{n-1}drd\om dt\\ \nonumber
=&\int_{0}^{T}\int_{2M}^{\infty}\int_{\mathbb{S}^{n-1}}ue^{-\la_0 t}\psi_{\la_0}\Big(\frac{r^{n-1}}{F}\big(\pa_{t}^{2}\eta_{T}^{2p'}-2\la_0\pa_{t}\eta_{T}^{2p'}\big)-Dr^{n-1}\pa_{t}\eta_{T}^{2p'}\Big)drd\om dt\ , \\ \label{725}
\les &\ T^{-1}\int_{T/3}^{T}\int_{2M}^{\infty}\int_{\mathbb{S}^{n-1}}|u|e^{-\la_0 t}\psi_{\la_0}\frac{r^{n-1}}{F}\eta_{T}^{2p'-2}drd\om dt
\end{align}
where $$C_{1}=\int_{2M}^{\infty}\int_{\mathbb{S}^{n-1}}r^{n-1}(\frac{u_{1}+\la_0 u_{0}}{F}+Du_{0})\psi_{\la_0}drd\om\ .$$
For simplicity, we denote 
$$A_1=\int_{0}^{T}\int_{2M}^{\infty}\int_{\mathbb{S}^{n-1}}|u|^{p}\eta_{T}^{2p'}\psi_{\la_0}e^{-\la_0 t}r^{n-1}drd\om dt\ .$$
By applying H\"older's inequality to \eqref{725}, we have that
\begin{align*}
&\int_{T/3}^{T}\int_{2M}^{\infty}\int_{\mathbb{S}^{n-1}}|u|e^{-\la_0 t}\psi_{\la_0}\frac{r^{n-1}}{F}\eta_{T}^{2p'-2}drd\om dt\\
\les &\ A_1^{\frac{1}{p}}
 \Big(\int_{\frac{T}{3}}^{T}\int_{2M}^{\infty}\int_{\mathbb{S}^{n-1}}F^{-p'}r^{n-1}e^{-\la_0 t}\psi_{\la_0}drd\om dt\Big)^{\frac{p-1}{p}}=A_1^{\frac{1}{p}}B_1^{\frac{p-1}{p}}\ .\\
\end{align*}
Recall that $dr=Fds$ and \eqref{srelation}, we get that
\begin{align*}
B_1\les&\int_{\frac{T}{3}}^{T}\Big(\int^{t+R}_{4M+e}\frac{r^{n-1}}{F^{p'-1}}e^{-\la_0 t}\psi_{\la_0}ds +\int_{-t-R}^{4M+e}\frac{r^{n-1}}{F^{p'-1}}e^{-\la_0 t}\psi_{\la_0}ds\Big) dt\\
\les &\int_{\frac{T}{3}}^{T}\Big(\int^{t+R}_{4M+e}s^{\frac{n-1}{2}}e^{-\la_0(t-s)}ds+\int_{-t-R}^{4M+e}e^{\frac{-s}{2M(p-1)}}e^{-\la_0 t}e^{s\la_0}ds\Big)dt\\
\les&\ T^{\frac{n+1}{2}}\ ,
\end{align*}
where we have used the fact that $\la_0>\frac{1}{2M(p-1)}$. Hence, we get that
$$\ep+A_1\ \les \ A_1^{\frac{1}{p}}T^{\frac{(n+1)(p-1)}{2p}-1}\ ,$$
which yields
\beeq
T\les \ \ep^{-\frac{2(p-1)}{n+1-(n-1)p}}, 
\eneq
for any $T\in[1, T_{\ep})$. We obtain the first part of lifespan estimate \eqref{d4}.

\subsection{$\frac{n}{n-1}<p<p_{S}(n)$}
In this subsection, we give the proof of the upper bound of lifespan estimate in the region $n/(n-1)<p<p_{S}(n)$. In the following, we take $T \in [\max\{8(4M+e), 36R\}, T_{\ep})$. 

Similar to the usual test function methods, we need to estimate the term 
$$L=\int_{0}^{T}\int_{2M}^{\infty}\int_{\mathbb{S}^{n-1}}|u|^{p}\eta_{T}^{2p'}(t) r^{n-1}drd\om dt \ .$$
For that purpose, we divide the integral into two parts by a smooth cut-off functions with respect to variable $s$ which was proposed by Lai-Zhou \cite{N-Y}
\begin{align*}
\al(s)=
\begin{cases}
\ \ 0,\ \ &-\infty\leq s\leq \frac{1}{8}\ ,\\
\ \ 1,\ \ &s\geq \frac{1}{4}\ .
\end{cases}
\end{align*}
Let $\al_{T}(s)=\al(s/T)$ and 
\begin{align*}
L_\infty=&\int_{0}^{T}\int_{2M}^{\infty}\int_{\mathbb{S}^{n-1}}|u|^{p}\eta_{T}^{2p'}(t)\al_{T}^{2p'}(s) r^{n-1}drd\om dt\ ,\\
L_{2M}=&\int_{0}^{T}\int_{2M}^{\infty}\int_{\mathbb{S}^{n-1}}|u|^{p}\eta_{T}^{2p'}(t)\big(1-\al_{T}^{2p'}(s)\big) r^{n-1}drd\om dt\ .
\end{align*}

Multiplying the equation in \eqref{h1} with $\eta_{T}^{2p'}(t)\al_{T}^{2p'}(s)\psi_{0} r^{n-1}$ where $\psi_0$ is in Lemma \ref{le6} and making integration by parts, then we come to
\begin{align*}
L_\infty\simeq &\int_{0}^{T}\int_{2M}^{\infty}\int_{\mathbb{S}^{n-1}}|u|^{p}\eta_{T}^{2p'}(t)\al_{T}^{2p'}(s)\psi_{0} r^{n-1}drd\om dt\\
=
&-\int_{0}^{T}\int_{2M}^{\infty}\int_{\mathbb{S}^{n-1}}\pa_{r}(r^{n-1}F)\psi_{0} \pa_{r}\al_{T}^{2p'}u\eta_{T}^{2p'}drd\om dt\\
&-\int_{0}^{T}\int_{2M}^{\infty}\int_{\mathbb{S}^{n-1}}r^{n-1}F(\psi_{0} \pa_{r}^{2}\al_{T}^{2p'}+2\pa_{r}\al_{T}^{2p'}\pa_{r}\psi_{0})u\eta_{T}^{2p'}drd\om dt\\
&+\int_{0}^{T}\int_{2M}^{\infty}\int_{\mathbb{S}^{n-1}}\frac{r^{n-1}}{F}(\pa_{t}^{2}\eta_{T}^{2p'})\al_{T}^{2p'}u\psi_{0} -ur^{n-1}D\psi_{0}(\pa_{t}\eta_{T}^{2p'})\al_{T}^{2p'}drd\om dt\ .
\end{align*}
Noticing that $|\pa_r(r^{n-1}F(r)) |\ \les \ r^{n-2}$, $|D|\ \les \ 1$ and 
$$|\pa^{2}_r\al_{T}^{2p'}| \les \ |\pa_r\al^{2p'}_T(s)| \les \ T^{-1}\al^{2p'-2}\ ,$$
$$ |\pa^{2}_t\eta_{T}^{2p'}| \les \ |\pa_t\eta^{2p'}_T(t)| \les \ T^{-1}\eta^{2p'-2}\ .$$
We get that  
$$L_{\infty}\ \les \ T^{-1}\int_{0}^{T}\int_{2M}^{\infty}\int_{\mathbb{S}^{n-1}}\frac{r^{n-1}}{F}\al^{2p'-2}\eta^{2p'-2} u\psi_{0} drd\om dt\ .$$
Then by H\"older's inequality, we have 
\begin{align*}
L_\infty 
\les &T^{-1}L_{\infty}^{\frac{1}{p}}\Big(\int_{\frac{T}{3}}^{T}\int_{4M+e}^{t+R}F^{1-p'}s^{n-1}dsdt\Big)^{\frac{1}{p'}}
\les\ T^{\frac{np-n-p-1}{p}} L_{\infty}^{\frac{1}{p}}\ ,
\end{align*}
which yields 
\beeq
\label{c1}
L_{\infty}\ \les \ T^{\frac{np-n-p-1}{p-1}}.
\eneq

Multiplying the equation in \eqref{h1} with $\eta_{T}^{2p'}(t)r^{n-1}$ and making integration by parts, then we come to
\begin{align*}
L=&\int_{0}^{T}\int_{2M}^{\infty}\int_{\mathbb{S}^{n-1}}|u|^{p}\eta_{T}^{2p'}(t)r^{n-1}drd\om dt\\
\leq&\int_{0}^{T}\int_{2M}^{\infty}\int_{\mathbb{S}^{n-1}}u\Big(\frac{r^{n-1}}{F}\pa_{t}^{2}\eta_{T}^{2p'}-Dr^{n-1}\pa_{t}\eta_{T}^{2p'}+r^{n-1}V\eta_{T}^{2p'}\Big)drd\om dt \\
\les &\int_{0}^{T}\int_{2M}^{\infty}\int_{\mathbb{S}^{n-1}}\frac{r^{n-1}}{F}|u| \eta_{T}^{2p'}drd\om dt\ .
\end{align*}
By applying H\"older's inequality, we get

\begin{align*}
L \ \les& \ L^{\frac{1}{p}}\Big(\int_{0}^{T}\int_{|s|\leq t+R}F^{1-p'}r^{n-1}dsdt\Big)^{\frac{p-1}{p}}\\
\les &\ L^{\frac{1}{p}}\Big(\int_{0}^{T}\int_{4M+e}^{t+R}s^{n-1}dsdt+\int_{0}^{T}\int_{-t-R}^{4M+e}e^{\frac{-s}{2M(p-1)}}dsdt\Big)^{\frac{p-1}{p}}\\
\les &\ L^{\frac{1}{p}}\Big(T^{n+1}+e^{\frac{T}{3M(p-1)}}\Big)^{\frac{p-1}{p}},
\end{align*}
which yields
\beeq
\label{c2}
L_{2M} \ \les\ L\ \les \ e^{\frac{T}{3M(p-1)}}\ .
\eneq
 
Next, as in \cite{N-Y}, we choose another test function 
$$\eta_{T}^{2p'}(t)\chi_{T}^{2p'}(t-s+2R)e^{-\la_0 t}\psi_{\la_0}r^{n-1}\ $$
where $\chi(\rho)$ is a smooth cut-off function satisfying 
\begin{align*}
\chi(\rho)=
\begin{cases}
\ 1\ ,&0\leq\rho\leq\frac{3}{4}\ ,\\
\ 0\ ,&\rho\geq\frac{5}{6}\ ,
\end{cases}
\end{align*}
and $\chi_{T}(\rho)=\chi(\rho/T)$. Thanks to \eqref{supp01}, it is easy to check that
\begin{align}
\nonumber
\begin{cases}
\supp\{\chi_{T}(2R-s)=1\} \ \cap \supp (u_{0},u_{1})=\supp (u_{0},u_{1}),\\ 
\supp\{\pa_{t}\chi_{T}(2R-s+t)\mid_{t=0}\} \cap \supp (u_{0},u_{1})=\emptyset \ . 
\end{cases}
\end{align}


\begin{figure}
\centering
\begin{tikzpicture}[scale=0.3]
\draw[->](-11,0)--(7,0) node[right] {s};
\draw[->](0,0)--(0,11) node[above] {t};
\draw(1,0)--(5,4);
\draw(-1,0)--(-7,6);
\draw(-4,0)--(2,6);
\draw(-8,0)--(0,8);
\draw[densely dotted](2,6)--(4,8);
\draw[densely dotted](0,8)--(2,10);
\draw[densely dotted](5,4)--(7,6);
\draw[densely dotted](-7,6)--(-8,7);
\filldraw[black!50!white] (-2.5,1.5)--(-4.5,3.5)--(0,8)--(2,6)--cycle;
\node[below] at (1,0) {\tiny $R$};
\node[below] at (-1,0) {\tiny $-R$};
\node[below] at (-4,0) {\tiny $\frac{-3T}{4}$$+$$2R$};
\node[below] at (-8,0) {\tiny $\frac{-5T}{6}$$+$$2R$};
\end{tikzpicture}
\caption{Support of the $\pa_{r}\chi_{T}(2R-s+t)$}
\label{tu1}
\end{figure}
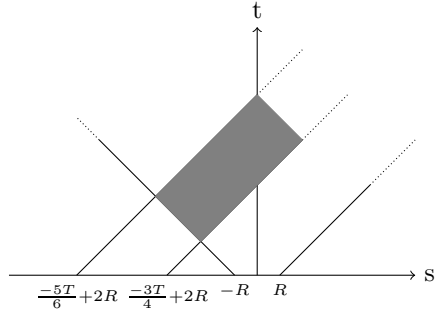
Multiplying the equation in \eqref{h1} with $\eta_{T}^{2p'}(t)\chi_{T}^{2p'}(t-s+2R)e^{-\la_0 t}\psi_{\la_0}r^{n-1}$ and making integration by parts, then we come to
\begin{align}\label{zhl}
&C_{2}\ep+\int_{0}^{T}\int_{2M}^{\infty}\int_{\mathbb{S}^{n-1}}|u|^{p}\eta_{T}^{2p'}(t)\chi_{T}^{2p'}(t-s+2R)e^{-\la_0 t}\psi_{\la_0}r^{n-1}drd\om dt\\ \nonumber
=&\int_{0}^{T}\int_{2M}^{\infty}\int_{\mathbb{S}^{n-1}}\frac{r^{n-1}}{F}u\psi_{\la_0}e^{-\la_0 t}\bigg\{\pa_{t}^{2}(\eta_{T}^{2p'}\chi_{T}^{2p'})-2\la_0\pa_{t}(\eta_{T}^{2p'}\chi_{T}^{2p'})-DF\pa_{t}(\eta_{T}^{2p'}\chi_{T}^{2p'})\\ \nonumber
& \ \ \ \ \ \ \ \ \ \ \ \ -\eta_{T}^{2p'}F^2\pa_{r}^{2}\chi_{T}^{2p'}-\eta_{T}^{2p'}F\frac{\pa_{r}(r^{n-1}F)}{r^{n-1}}\pa_{r}\chi_{T}^{2p'}-2\eta_{T}^{2p'}F^2\frac{\pa_{r}\psi_{\la_0}}{\psi_{\la_0}}\pa_{r}\chi_{T}^{2p'}\bigg\}drd\om dt\\  \nonumber
=&\int_{T/3}^{T}\int_{2M}^{\infty}\int_{\mathbb{S}^{n-1}}\frac{r^{n-1}}{F}u\psi_{\la_0}e^{-\la_0 t}K(t, r)dr d\om dt:= A_2\ ,
\end{align}
where $$C_{2}=\int_{2M}^{\infty}\int_{\mathbb{S}^{n-1}}r^{n-1}\psi_{\la_0}(\frac{u_{1}+\la_0 u_{0}}{F}+Du_{0})d\om dr\ .$$
To estimate the term $A_2$, we first note that the integral interval of time is contained in $(T/3, T)$. In fact, by the definition of $\chi$, it is easy to check that the support of $\pa_{r}\chi_{T}(t-s+2R)$ is
\begin{equation}\label{73102}
\frac{3T}{4}\leq t-s+2R\leq \frac{5T}{6}\ .
\end{equation}
 Recalling that the support of solution $u$ satisfying \eqref{fps}, then we can get that $t\geq T/3$ by combing \eqref{73102} and \eqref{fps}. Besides, by \eqref{b4}, we see that $(\pa_{r}\psi_{\la_0})/(\psi_{\la_0}) \les F^{-1}$, so we have 
 $$K(t, r)\  \les\ T^{-1}\eta_{T}^{2p'-2}(t)\chi_{T}^{2p'-2}\ ,$$
 due to the fact that $0\leq F\leq1$. 
Next we divide $A_2$ into two parts
\begin{align*}
A_2\ \leq\  &\int_{0}^{T}\int_{2M}^{\infty}\int_{\mathbb{S}^{n-1}}|u|\frac{r^{n-1}}{F}\psi_{\la_0}e^{-\la_0 t} K(t, r)|\al_{T}^{2p'}(s)|drd\om dt\\
+&\int_{0}^{T}\int_{2M}^{\infty}\int_{\mathbb{S}^{n-1}}|u|\frac{r^{n-1}}{F}\psi_{\la_0}e^{-\la_0 t} K(t, r)|1-\al_{T}^{2p'}(s)|drd\om dt\\
=&: A^{1}_2+A_2^2\ .
\end{align*} 
By applying H\"older's inequality to $A^{1}_2$ and $A^{2}_2$, we get that

\begin{align*}
|A_{2}^{1}|\ \les\ & T^{-1}L_{\infty}^{\frac{1}{p}}\Big(\int_{\frac{T}{3}}^{T}\int_{2M}^{\infty}e^{-\la_0 p'(t-s)}F^{-p'}r^{n-1-\frac{(n-1)p'}{2}}drdt\Big)^{\frac{p-1}{p}}\\
\les&T^{-1}L_{\infty}^{\frac{1}{p}}\Big(\int_{\frac{T}{3}}^{T}\int_{4M+e}^{t+R}e^{-\la_0 p'(t-s)}F^{1-p'}s^{n-1-\frac{(n-1)p'}{2}}dsdt\Big)^{\frac{p-1}{p}}\\
\les &L_{\infty}^{\frac{1}{p}}T^{(n-\frac{(n-1)p'}{2})\frac{1}{p'}-1}\ ,
\end{align*}
\begin{align*}
|A_{2}^{2}|\ \les \ &T^{-1}L_{2M}^{\frac{1}{p}}\Bigg(\int_{\frac{T}{3}}^{T}\int_{-t-R}^{4M+e}e^{-\la_0 p't}e^{sp'\la_0}F^{1-p'}r^{n-1-\frac{(n-1)p'}{2}}dsdt\\
&\ \ \ \ \ \ \ \ \ \ \ \ +\int_{\frac{T}{3}}^{T}\int_{4M+e}^{\frac{T}{4}}e^{-\la_0 p'(t-s)}F^{1-p'}r^{n-1-\frac{(n-1)p'}{2}}dsdt\Bigg)^{\frac{p-1}{p}}\\
\ \les\ &T^{-1}L_{2M}^{\frac{1}{p}}\Bigg(\int_{\frac{T}{3}}^{T}T^{n-1-\frac{(n-1)p'}{2}}e^{-\la_0 p' t}dt\int_{4M+e}^{\frac{T}{4}}e^{\la_0 p's}ds\\
&\ \ \ \ \ \ \ \ \ \ \ \ \ \ +\int_{\frac{T}{3}}^{T}e^{-\la_0 p't}dt\int_{-t-R}^{4M+e}e^{s\la_0 p'-\frac{s}{2M(p-1)}}ds\Bigg)^{\frac{p-1}{p}}\\
\ \les\ &T^{-1}L_{2M}^{\frac{1}{p}}\Big(T^{n-\frac{(n-1)p'}{2}}e^{\frac{-\la_0 Tp'}{12}}+T^{2}e^{\frac{-\la_0 p'T}{3}}\Big)^{\frac{p-1}{p}}\ .
\end{align*}
By combing \eqref{c1} \eqref{c2} and \eqref{zhl}, we get that
\begin{align}\nonumber  
C_{2}\ep\ \leq \ A^2_2\ \les\ & L_{\infty}^{\frac{1}{p}}T^{-1+(n-\frac{(n-1)p'}{2})\frac{1}{p'}}+L_{2M}^{\frac{1}{p}}T^{-1}(T^{n-\frac{(n-1)p'}{2}}e^{\frac{-\la_0 p'T}{12}}+T^{2}e^{\frac{-\la_0 p'T}{3}})^{\frac{1}{p'}}\\  \label{c3}
\les\ &\ L_{\infty}^{\frac{1}{p}}T^{-1+(n-\frac{(n-1)p'}{2})\frac{1}{p'}}+T^{\frac{(n-1)p^{2}-3p(n+1)+2n}{2p(p-1)}}e^{T(\frac{1}{3M(p-1)}-\frac{\la_0}{12})}\\
\les\ &T^{\frac{(n-1)p^{2}-(n+1)p-2}{2p(p-1)}}, \nonumber
\end{align}
 for any $T \in [\max\{8(4M+e), 36R\}, T_{\ep})$. Hence we obtain the second lifespan estimate of \eqref{d4}.

\subsection{$p=p_{S}(n)$}

In this subsection, we give the proof of the upper bound lifespan estimates when $p=p_S(n)$. We choose $b_{a}\al_{T}^{2p'}(s)\eta_{T}^{2p'}(t)r^{n-1}$ as test function where $b_a$ is in \eqref{bade} and $a=\frac{n-1}{2}-\frac{1}{p}$, $T \in [\max\{8(4M+e), 36R\}, T_{\ep})$.  

Multiplying the equation in \eqref{h1} with $b_{a}\al_{T}^{2p'}(s)\eta_{T}^{2p'}(t)r^{n-1}$ and making integration by parts, then we come to
\begin{align}\nonumber
X(T)=&\int_{T/16}^{T}\int_{2M}^{\infty}\int_{\mathbb{S}^{n-1}}|u|^{p}b_{a}\eta^{2p'}_{T}(t)\al_{T}^{2p'}(s)r^{n-1}drd\om dt\\ \nonumber
=&\int_{T/16}^{T}\int_{2M}^{\infty}\int_{\mathbb{S}^{n-1}}\frac{r^{n-1}}{F}u b_a\bigg(\al_{T}^{2p'}\pa_{t}^{2}\eta_{T}^{2p'}-2\al_{T}^{2p'}\pa_{t}\eta_{T}^{2p'}\frac{b_{a+1}}{b_a}-\eta_{T}^{2p'}\frac{F\pa_{r}(r^{n-1}F)}{r^{n-1}}\pa_{r}\al_{T}^{2p'}\\ \nonumber
&\ \ \ \ \ \ \ \ \ \ \ \ \ \ \ \ \ \ \ -F^2\eta_{T}^{2p'}\pa_{r}^{2}(\al_{T}^{2p'})-2F^2\eta_{T}^{2p'}\pa_{r}\al_{T}^{2p'}\frac{\pa_{r}b_a}{b_a}-DF\pa_{t}\eta_{T}^{2p'}\al_{T}^{2p'}
\bigg)drd\om dt\\ \label{khbd}
=&\int_{T/16}^{T}\int_{2M}^{\infty}\int_{\mathbb{S}^{n-1}}\frac{r^{n-1}}{F}u b_a U(t, r)drd\om dt\ .
\end{align}
Note that the support of $u\al^{2p'}_{T}(s)$ with respect to $t$ is contained in $(T/16, T)$. In fact, by the support of $\al_{T}^{2p'}(s)$ and \eqref{fps} we have
$$\frac{T}{8}\leq s\leq t+R\ ,$$
thus $t\geq T/8-R\geq T/16$ since $T\geq 16R$. To estimate $U(t, r)$, we see that 
$$D(r)\ \les\  r^{-1}\  \les \ T^{-1}\ , \ \frac{\pa_{r}(r^{n-1}F)}{r^{n-1}}\ \les\ r^{-1} \ \les \ T^{-1}\ , \ T/8\leq s\leq T+R\ .$$
Besides, by Lemma \ref{leba}, we see that $|\pa_r b_a| \les b_{a+1}$ and 
$$\frac{b_{a+1}}{b_a}\gtrsim \frac{(t+R)^{-(a+1)}}{(t+R)^{-a}}\simeq T^{-1}\ , t\in(T/16, T)\ . $$
So we obtain that 
\beeq\label{khbd1}
|U(t, r)|\ \les\ T^{-1}\al_{T}^{2p'-2}\eta_{T}^{2p'-2}\frac{b_{a+1}}{b_a}\ .
\eneq
By applying H\"older's inequality to \eqref{khbd} and using \eqref{khbd1}, we have \begin{align*}
X(T)\ \les\ & T^{-1}X(T)^{\frac{1}{p}}
\Big(\int_{\frac{T}{16}}^{T}\int_{4M+e}^{t+R}b_{a}^{\frac{-1}{p-1}}b_{a+1}^{\frac{p}{p-1}}s^{n-1}F^{1-p'}dsdt\Big)^{\frac{p-1}{p}}\\
\les\ &T^{-1} X(T)^{\frac{1}{p}}\Big(\int_{\frac{T}{16}}^{T}\int_{4M+e}^{t+R}(t+R+s)^{\frac{n-1}{2}-\frac{1}{p(p-1)}}(t+R+1-s)^{-1}dsdt\Big)^{\frac{p-1}{p}}\\
\les\ &X(T)^{\frac{1}{p}}(\ln T)^{\frac{p-1}{p}},
\end{align*}
which yields
\begin{align}\label{khbd3}
X(T)\ \les \ (\ln T)^{1/p'}X(T)^{1/p}\ .
\end{align}
As in \cite{L-L-W-T}, we introduce 
\beeq\nonumber
Y(N)=\int_{16}^{N}X(T)T^{-1}dT\ ,\ N>16\ , 
\eneq
then it is to see that 
$$Y'(N)=N^{-1}X(N)\ .$$
Furthermore, by doing direct calculations we have
\begin{align}\nonumber
Y(N)\les&\int_{1}^{N}\int_{2M}^{\infty}\int_{\mathbb{S}^{n-1}}|u|^{p}b_a\al_{N}^{2p'}r^{n-1}\int_{t}^{\min (N,16t)}(\eta_{T}(t))^{2p'}T^{-1}dTdtdrd\om\\  \nonumber\les&\int_{1}^{N}\int_{2M}^{\infty}\int_{\mathbb{S}^{n-1}}|u|^{p}b_a\al_{N}^{2p'}r^{n-1}\int^{1}_{\max (\frac{1}{16},\frac{t}{N})}(\eta(\tau))^{2p'}\tau^{-1}d\tau dtdrd\om\\ \label{khbd4}
\les&\int_{1}^{N}\int_{2M}^{\infty}\int_{\mathbb{S}^{n-1}}|u|^{p}b_a\al_{N}^{2p'}r^{n-1}(\eta_N(t))^{2p'}\int_{\frac{1}{16}}^{1}\tau^{-1}d\tau dtdrd\om\\ \nonumber\les&\ln16\int_{0}^{N}\int_{2M}^{\infty}\int_{\mathbb{S}^{n-1}}|u|^{p}b_a\al_{N}^{2p'}r^{n-1}(\eta_N(t))^{2p'}dtdrd\om\\ \nonumber
\les& X(N)
\end{align}
Thus, by combing \eqref{khbd3}-\eqref{khbd4}, we have
\beeq
\label{008}
Y(T)\les(\ln T)^{1/p'}(TY'(T))^{1/p}.
\eneq
In addition, by \eqref{c3}, we can get the lower bound of $L_\infty$
\beeq\nonumber
L_\infty \gtrsim \ep^{p}T^{n-\frac{(n-1)p}{2}} \ ,
\eneq
then by the lower bound of $b_a$, we can get
\begin{align} \nonumber 
TY'(T)=&\int_{\frac{T}{16}}^{T}\int_{2M}^{\infty}\int_{\mathbb{S}^{n-1}}|u|^{p}b_{a}\eta^{2p'}_{T}\al_{T}^{2p'}(s)r^{n-1}drd\om dt\\ \nonumber
\gtrsim &\ L_{\infty}T^{-\frac{n-1}{2}+\frac{1}{p}}\\ \label{c7}
\gtrsim & \ T^{n-\frac{(n-1)p}{2}-\frac{n-1}{2}+\frac{1}{p}}\ep^{p}=\ep^{p}\ ,
\end{align}
where we used the fact 
$$n-\frac{(n-1)p}{2}=\frac{n-1}{2}-\frac{1}{p}$$
when $p=p_{S}(n)$.
By applying Lemma \ref{le11} to \eqref{c7} and \eqref{008}, we obtain the last lifespan estimate of Theorem \ref{thm1}\ . 
\section{Proof of Theorems \ref{thm2}}
In this section, we give the proof of Theorem \ref{thm2}. We shall transform \eqref{h1} into an one dimensional wave equation by using Regge-Wheeler coordinate \eqref{wg}.

\subsection{Transformation}We shall consider the radial solutions $u=u(t,s)$, 
then \eqref{h1} can be written as 
\beeq
\nonumber
\pa^{2}_{t}u-\pa^{2}_{s}u-\frac{(n-1)F}{r(s)}\pa_{s}u+D(r(s))Fu_{t}+V(r(s))Fu=F|u_t|^{p}\ .
\eneq
Let $u(t,s)=v(t,s)r(s)^{-\frac{n-1}{2}}$, then $v(t, s)$ satisfying 
\begin{align}
\label{h8}
\begin{cases}
\pa_{t}^{2}v-\pa_{s}^{2}v+Q(r(s))v+F(s)D((s))v_{t}=FJ(s)^{p-1}|v_t|^{p}\ , \\
v(0,s)=\ep f(s),\ \ \ v_{t}(0,s)=\ep g(s)\ ,
\end{cases}
\end{align}
where 
\beeq
\label{R1}
J=r(s)^{-\frac{n-1}{2}},\ Q=F(r(s))V(r(s))+\frac{(n^{2}-4n+3)F^{2}}{4r(s)^{2}}+\frac{(n-1)MF}{r(s)^{3}}\ .
\eneq
In addition, by the relation between $s$ and $r$ \eqref{srelation}, we have
\begin{align}
FJ(s)^{p-1}\simeq
\begin{cases}
s^{\frac{(1-p)(n-1)}{2}}\ ,\ \ s\geqslant4M+e\ ,\\
e^{\frac{s}{2M}}\ ,\ \ \ \ \ \ \ \ \ \ s\leqslant4M+e\ .
\end{cases}
\end{align}
\subsection{Test function}
 \begin{lem}
\label{le4}
Suppose $D, V\in C(2M, \infty)$ and satisfying  \eqref{khbd00}. Then for any $\la>0$, the following equation 
\beeq
\pa_{s}^{2}\phi_{\la}=Q\phi_{\la}+\la DF\phi_{\la}+\la^{2}\phi_{\la}\ ,
\eneq
admits a $C^{2}$ solution satisfying
\begin{align}
\label{b1}
\phi_{\la}\simeq e^{\la s},
\end{align}
\end{lem}
\begin{prf}
When $s>0$, we consider the following ordinary differential equation
\begin{align}
\begin{cases}
\phi''_{\la}-\la^{2}\phi_{\la}=(Q+\la DF)\phi_{\la},\\
\phi_{\la}(0)=1,\ \ \phi'(0)=\la.
\end{cases}
\end{align}
Let $Y=(\phi_{\la}\ \phi'_{\la})^{T}$, then we can get that
\begin{align*}
Y'=(A+B(s))Y,\ \  A=\begin{pmatrix} 0 & 1\\ \la^{2} & 0 \end{pmatrix},\ \ B=\begin{pmatrix} 0 & 0\\ Q+\la DF & 0\end{pmatrix},
\end{align*}
where $B\in L^{1}[0,\infty)$. By applying Levinson theorem (see, e.g., [\cite{E-A-C} chapter 3, Theorem 8.1]) to the system. $B(s)$ is negligible perturbation term. Then there exists $s_{0} \in [0,\infty)$ such that
\begin{align*}
Y\simeq c_{3}\begin{pmatrix} 1\\ \la \end{pmatrix}e^{\la s}+c_{4}\begin{pmatrix} 1\\ -\la \end{pmatrix}e^{-\la s}, s\geq s_{0},
\end{align*}
where $c_{3}$ and $c_{4}$ are constant that can be chosen. Concerning the initial data, we can choose $c_{4}=0$, that is to say $\phi_{\la}\sim e^{\la s}$.
For the remaining case, set $X(s)=\phi_{\la}(-s)$, we consider the following ordinary differential equation problem
\begin{align*}
\begin{cases}
X''-\la^{2}X=(Q+\la DF)X\ ,\\
X(0)=1\ ,\ \ \ X'(0)=-\la, s>0\ ,
\end{cases}
\end{align*}
Samilarly, we can obtain $X(s)\simeq e^{-\la s}$, that is $\phi_{\la}\sim e^{\la s}$.
\end{prf}

With the help of Lemma \ref{le4}, it is easy to see that $\Phi(t,s)=e^{-\la t}\phi_{\la}$ is the solution of Linear wave equation
$$\pa^{2}_{t}u-\pa_{s}^{2}u-DF\pa_{t}u+Qu=0\ ,$$
for any $\la>0$.

\subsection{Lifespan estimates for $1<p<p_G(n)$} Let $\la=\la_{0}>\frac{4}{Mp(p-1)}$ to be taken and $\Phi=\phi_{\la_0}e^{-\la_0 t}$, where $\phi_{\la_0}$ is in Lemma \ref{le4}. By multiplying the equation in \eqref{h8} with $-\pa_{t}\big(\eta^{2p^{\prime}}_{T}(t)\Phi\big)$ and making integration by parts, then we have 
\begin{align}\label{khbd7}
&C_{3}\ep+\int_{0}^{T}\int_{\R}F|v_{t}|^{p}J(s)^{p-1}\Big(\eta_{T}^{2p'}\Phi \la_0-2p'\eta_{T}^{2p'-1}(\pa_{t}\eta_{T})\Phi\Big) dsdt\\ \nonumber
= &\int_{0}^{T}\int_{\R}v_t\Big((\pa_{t}^{2}\eta_{T}^{2p'})\Phi+2(\pa_{t}\eta_{T}^{2p'})(\pa_{t}\Phi)-(\pa_{t}\eta_{T}^{2p'})\Phi DF\Big)dsdt  \\  \nonumber
\les&\  T^{-1}\int_{\frac{T}{3}}^{T}\int_{\R}F|v_t|\Phi(t,s) \eta_{T}^{2p'-2}(t)dsdt\ ,
\end{align}
where $C_3=\int_{\R}\phi_{\la_0}\Big(g+(\la_0 DF+\la^2_{0}-Q)f\Big) ds >0$ due to the assumption \eqref{chuzhi2}. 
 We denote 
$$X(T)=\int^{T}_{T/3}\int_{\R}|v_{t}|^{p}\Phi \eta_{T}^{2p'}FJ^{p-1}dsdt \ , \ W(T)=\int^{T}_{0}\int_{\R}|v_{t}|^{p}\Phi \eta_{T}^{2p'}FJ^{p-1}dsdt\ ,$$
then it is easy to see $X(T)\leq W(T)$. By \eqref{khbd7}, we have 
\beeq\label{h09}
C_3\ep+W(T) \ \les \ T^{-1}\int_{\frac{T}{3}}^{T}\int_{\R}F|v_t|\Phi(t,s) \eta_{T}^{2p'-2}(t)dsdt\ ,
\eneq
due to $\eta'(t)\leq 0$. By applying H\"older inequality to \eqref{h09}, we have that
 \begin{align}\nonumber
 &C_{3}\ep+W(T)\\ \nonumber
 \ \les \ &T^{-1}X(T)^{\frac{1}{p}}\Big(\int_{\frac{T}{3}}^{T}\int_{\R}\Phi FJ^{-1}dsdt\Big)^{\frac{p-1}{p}} \\ \nonumber
\les & T^{-1}X(T)^{\frac{1}{p}}\Bigg(\int_{\frac{T}{3}}^{T}\Big(\int_{4M+e}^{t+R}s^{\frac{n-1}{2}}e^{-\la_0(t-s)}ds+\int_{-t-R}^{4M+e}e^{-\la_{0}(t-s)}e^{\frac{s}{2M}}ds\Big)dt\Bigg)^{\frac{p-1}{p}}\\  \label{khbd9}
\les&  \ X(T)^{\frac{1}{p}}T^{\frac{n+1}{2p'}-1} \ .
\end{align}
By applying Young's inequality to \eqref{khbd9}, we get that
$$C_{3}\ep+W(T)\leq  \frac{X(T)}{p}+C T^{\frac{n+1}{2}-p'}\leq \frac{W(T)}{p}+C T^{\frac{n+1}{2}-p'}$$
for some $C>0$, which yields 
\begin{align}
T_{\ep}\ \les \ \ep^{-(\frac{1}{p-1}-\frac{n-1}{2})^{-1}} \ .
\end{align}

\subsection{Lifespan estimate for $p=p_{G}(n)$}
We define a functional 
$$Y(N)=\int_{3}^{N}X(T)T^{-1}dT\ , N\geq 3\ .$$
Then by the same procedure in \eqref{khbd4}, we have 
$$Y(T) \ \les\  W(T)\ .$$
According to \eqref{khbd9}, we see that 
\beeq
\label{8}
\big(C_{3}\ep+Y(T))\big)^{p} \lesssim TY'(T)\ ,
\eneq
where we have used the fact 
$$\frac{(n+1)p}{2p'}-p=0\ , \ p=p_G(n)\ .$$
Let $Z(T)=C_{3}\ep+Y(T)$, then $Z(3)=C_{3}\ep$ and 
$$TZ'(T)\ \gt \ Z(T)^{p}\ , Z(3)=C_{3}\ep\ .$$
This is a typical ODE inequality which will blow up in finite time and we shall have
$$T_{\ep} \ \les\ \exp \big(\ep^{-(p-1)}\big)\ .$$

\bibliographystyle{plain1}

\end{document}